\documentstyle[12pt]{amsart}
\input epsf

\newtheorem{theorem}{Theorem}[section]
\newtheorem{definition}{Definition}[section]

\newtheorem{lemma}{Lemma}[section]

\newtheorem{corollary}{Corollary}[section]
\newenvironment{proof}{\noindent {\bf Proof:}}{$\Box$}

\newenvironment{ppf}{\noindent {\bf Proof of Theorems 1.1 and 1.2:}}{$\Box$}

\numberwithin{equation}{section}
\numberwithin{figure}{section}

\newcommand\ZHS{$\Bbb Z$-homology 3-sphere }
\newcommand\ZHSs{$\Bbb Z$-homology 3-spheres }

\begin{document}

\title[Finite type 3-manifold invariants]{On Vassiliev knot invariants induced
\\ from finite type
	3-manifold invariants}
\author{Matt Greenwood}
\address{Department of Mathematics, Columbia University, New York, NY 10027}
\email{matt@@math.columbia.edu}
\author{Xiao-Song Lin}
\address{Department of Mathematics, Columbia University, New York, NY 10027}
\email{xl@@math.columbia.edu}
\date{May, 1995}
\thanks{The second author is supported in part by NSF.}

\maketitle

\begin{abstract} We prove that the knot invariant induced by a \ZHS invariant
of
order $\leq k$ in Ohtsuki's sense, where $k\geq4$, is of order $\leq k-2$. The
method developed in our computation shows that there is no \ZHS invariant of
order 5.
\end{abstract}

\section{Introduction}

The theory of finite type knot invariants is by now well developed
(see [V], [Gu], [B-L], [K] and [BN]). It provides an unified approach
to various quantum knot invariants (Jones, HOMFLY and Kauffman
polynomials, etc.). Such an analogous unified approach to quantum
3-manifold invariants is not available at present.

In [O1], Ohtsuki reorganized quantum $SU(2)$ invariants $\tau_r(M)$
for a \ZHS $M$ into a formal power series
$\tau(M)=\sum_{n=0}^\infty\lambda_n(M)(t-1)^n$. In order to find a
characterization of the invariants $\lambda_n(M)$, Ohtsuki introduced
the notion of finite type \ZHS invariants. Precisely, let $L$ be an
algebraically split link with framings $\pm 1$ on each component in a
\ZHS $M$ and we will denote by $M_{L}$ the {\ZHS} obtained by surgery
on $L$. Let $\lambda$ be a \ZHS invariant, then $\lambda$ is said to
be a {\it finite type invariant} of order $\leq k$ if the following
sum vanishes for all choices of $M$ and $L$ such that $\#L>k$:

\begin{equation}
\sum_{L' \subset L}(-1)^{\#L'}\lambda(M_{L'}) = 0.
\end{equation}
Here $\#L$ is the number of components of the link $L$.

Let ${\cal O}_k$ be the vector space of all finite type invariants of
\ZHSs of order $\leq k$. The main result of [O2] is that
$\text{dim}({\cal O}_k)<\infty$. Moreover, it is proved that the \lq\lq
$k$-th difference" of a finite type \ZHS invariant of order $k$ is
determined by its \lq\lq weight" on the set of graphs of $k$ edges
such that (1) the valence at each vertex is either 1 or 3 and (2)
every edge connects two distinct vertices.

In [G], among other things, Garoufalidis proved that ${\cal O}_4={\cal
O}_3$. He also observed that every $\lambda\in{\cal O}_k$ determines a
finite type knot invariant $v_\lambda$ of order $\leq k-1$ by
\begin{equation}v_\lambda(K)=\lambda(S^3_K)
\end{equation}
for every knot $K$ in $S^3$ with $+1$
framing. A way of calculating the weight system of $v_\lambda$ was
also described in [G]. Garoufalidis suggested that it might be
interesting to study these induced knot invariants, although it was
also observed that the presence of Kirby moves may restrict the
induced knot invariants considerably.

In this paper, we examine
closely Ohtsuki's proof of the finiteness of $\text{dim}({\cal
O}_k)$. This leads to one of our main results showing that the order of a
induced knot invariant is, in most cases, actually one less than what it was
supposed to be {\it a priori}.

\begin{theorem} Let $v_\lambda$ be the knot invariant induced by
$\lambda\in{\cal O}_k$ with $k\geq 4$. Then $v_\lambda$ is of order at most
$k-2$.
\end{theorem}

Our method also leads to the determination of \ZHS invariants of order $\leq
5$.

\begin{theorem} ${\cal O}_5={\cal O}_4$.
\end{theorem}

In [R1,2], Rozansky has shown, using
physical predictions of the asymptotic behavior of Witten's Chern-Simons
path integral over trivial connections, that  Ohtsuki's invariant
$\lambda_n$ is of
order $3n$ and he further conjectured that ${\cal O}_{3n}={\cal O}_{3n+1}={\cal
O}_{3n+2}$. The results of this paper certainly point in this direction.

Let $\lambda_{\text{C}}$ be the Casson invariant and 1 be the constant
invariant whose value on each \ZHS is 1. It is known that
$\{1,\lambda_{\text{C}}\}$ is a basis for ${\cal O}_3$ (see [O2] and
[G]). Let $v_{\text{C}}$ be the knot invariant induced by
$\lambda_{\text{C}}$. It is known that $v_{\text{C}}$ is of order 2. To
illustrate Theorem 1.1, we may consider the order 6 \ZHS invariant
$\lambda_{\text{C}}^2$ whose induced knot invariant $v_{\text{C}}^2$ is of
order
4. Probably also motivated by this example, Garoufalidis asked whether the knot
invariant induced by a \ZHS invariant of order $\leq 3n$ is actually of order
$\leq 2n$. The answer to this question remains open.

\section{Definitions and basic properties.}

In this section we will first recall a few definitions and lemmas of
Ohtsuki [O2]. Then we will give a proof of Ohtsuki's finiteness theorem
using the notion of surgery classification of links in $S^3$. Our
proof here will have the same essential ingredients as in Ohtsuki's
proof. But we do think that our exposition is clearer than that
of Ohtsuki's. We single out a key step in our argument, the proof of
Theorem 2.2, in the next section since the same method will be used in
the proof of Theorem 1.1 in the last section.

An {\it algebraically split link} (ASL) $L$ in a \ZHS $M$ is
a link with unoriented and unordered
components such that all of its linking numbers are zero.
Henceforth, all 3-manifolds will be \ZHSs and all links
will be $\pm1$-framed ASL's.

Let $\cal I$ be the set of homeomorphism classes of $\Bbb Z$-homology
3-spheres. Let $\cal V$ be the vector space over $\Bbb C$ with $\cal
I$ as a basis. Define an element $(M,L)$ in $\cal V$ by
$$(M, L) =\sum_{L' \subset L}(-1)^{\#L'}
M_{L'}$$
where $L'$ runs through all sublinks (including the empty sublink) of $L$.
Let ${\cal V}_k$ be the subspace of $\cal V$ spanned by $(
M,L)$ for all \ZHSs $M$ and all $\pm1$-framed links $L$ such
that $\#L = k+1$. The space $\cal V$ has a natural
stratification:
$$\cal V={\cal V}_{-1}\supset{\cal V}_0\supset{\cal
V}_1\supset\cdots\supset{\cal V}_{k-1}\supset{\cal V}_k\supset\cdots.$$

For $A,B\in{\cal V}_{k-1}$, we denote $A\sim B$ if $A-B\in{\cal
V}_k$. We have the following basic lemmas from [O2].

\begin{lemma} For any $(M,L)\in{\cal V}_{k-1}$, there is a link $J$ in $S^3$
with
$\#J=k$ such that $(M,L)\sim(S^3,J)$.
\end{lemma}

\begin{lemma} $(S^3,L)\sim(S^3,J)$ if $J$ can be obtained from
$L$ by either of the local moves shown in Figure 2.1 where $L_i$ and $
L_j$ may be different components of $L$. In the second local move, we may
orient components of $L$ in any favorable way.
\end{lemma}

\begin{lemma} Let $L$ and $J$ be the same link in $S^3$ with different
$\pm1$ framings. Let $s(L)$ be the product of framings of $L$. Then
$s(L)(S^3,L)\sim s(J)(S^3,J)$.
\end{lemma}

\begin{figure}
\caption{\protect\footnotesize\protect\sf Local moves on ASL's.}
\end{figure}

\bigskip

{}From these lemmas, we may think of ${\cal V}_{k-1}/{\cal V}_k$ as spanned by
equivalence classes of ASL's in $S^3$ with $k$ components under the
equivalence relations depicted in Figure 2.1, where an ASL $L$ is
regarded as an element $(S^3,L)\in{\cal V}_{k-1}$ with the understanding
that the framings of $L$ are all $+1$.  Now it seems to be the right
moment to recall the notion of surgery classification of links in
$S^3$ here.

\begin{definition} A {\em surgery modification} on an ASL $L$ in $S^3$
is first to have a disk $B$ in $S^3$ intersecting $L$ only in its
interior with zero algebraic intersection number. Then, perform a
$\pm1$ surgery on $\partial B$ which will change $L$ to another ASL
$J$ in $S^3$. Two ASL's in $S^3$ are {\em surgery equivalent} if one
of them can be changed to another by a finite sequence of surgery
modifications.
\end{definition}

Notice that Figure 2.1 also serves as an illustration of a surgery
modification. Thus, we have

\begin{lemma} If $L$ and $J$ are surgery equivalent ASL's with $k$ components,
the $(S^3,L)\sim(S^3,J)$.
\end{lemma}

To arrive at the finiteness of $\text{dim}({\cal V}_{k-1}/{\cal V}_k)$, we
recall the following theorem of Levine [L].

\begin{theorem} Surgery equivalence classes of ASL's with oriented
and ordered components are classified by Milnor's triple $\mu$-invariants.
\end{theorem}

We may now think of a surgery equivalence class of ASL's with $k$
components as a collection of integers $\{\mu(i_1i_2i_3)\}$ where
$i_1,i_2,i_3$ are distinct indices among $1,\dots,k$. They satisfy the
following relations:
$$\left\{\begin{aligned}
&\mu(i_1i_2i_3)=\mu(i_2i_3i_1)=\mu(i_3i_1i_2) \\
&\mu(i_1i_2i_3)=-\mu(i_2i_1i_3)
\end{aligned}
\right. $$
Take the
advantage of being in the vector space
${\cal V}_{k-1}/{\cal V}_k$, each collection of integers
$\{\mu(i_1i_2i_3)\}$ will be further simplified in the following theorem. This
theorem will imply the finiteness of
$\text{dim}({\cal V}_{k-1}/{\cal V}_k)$.

\begin{theorem} In ${\cal V}_{k-1}/{\cal V}_k$, an ASL representing
$\{\mu(i_1i_2i_3)\}$ can be expressed as a linear combination of ASL's whose
collections of triple $\mu$-invariants satisfy the following conditions:
\begin{enumerate}
\item each $\mu(i_1i_2i_3)$ is either $0$ or $1$ for certain fixed
cyclic order of indices;
\item each index $i$ appears in at most two non-zero $\mu(i_1i_2i_3)$'s.
\end{enumerate}
\end{theorem}

Now collections $\{\mu(i_1i_2i_3)\}$ of $0$'s and $1$'s, where
$i_1,i_2,i_3$ are distinct indices among $1,\dots,k$ with a fixed
cyclic ordering are in 1-1 correspondence with graphs with $k$ ordered
edges such that
\begin{enumerate}
\item each vertex is of valence either 1 or 3;
\item each edge connects two distinct vertex.
\end{enumerate}
Furthermore, we may forget the ordering of edges of such a graph since
the components of ASL's in ${\cal V}_{k-1}/{\cal V}_k$ need not to be
ordered. Thus, we see that ${\cal V}_{k-1}/{\cal V}_k$ is spanned by
ASL's in 1-1 correspondence with graphs having $k$ edges and
satisfying the properties (1) and (2) above. In particular, we have
the following corollary.

\begin{corollary} $\text{dim}({\cal V}_{k-1}/{\cal V}_k)<\infty$.
\end{corollary}

\section{The proof of Theorem 2.2}

We start with an ASL $L$ in $S^3$ representing
$\{\mu(i_1i_2i_3)\}$. We may assume that each component of $L$ is
unknotted. Notice that modifying one component of $L$ will not effect
the triple $\mu$-invariants not involving the index of that
component. This allows us to \lq\lq localize" a component of $L$ in
the following way.

Pick a component, say $L_1$, of $L$. It bounds a disk $B$. Other
components of $L$ will intersect $B$ in its interior only and the
algebraic intersection numbers are all zero. Thus, we may pair the
intersection points so that each pair consisting of two points lying
on the same component of $L$, one with positive intersection number
and the other with negative intersection number. On the disk $B$, we
draw small circles $w_1,w_2,\dots,w_n$ such that
\begin{enumerate}
\item each of them surrounds a pair of intersection points;
\item they are all disjoint; and
\item they are not nested with each other.
\end{enumerate}
We may further assume that each small circle $w_i$ represents a simple
commutator $[m_{i_1},m_{i_2}]$, where $m_{i_1}$ and $m_{i_2}$ are
meridians of the $i_1$-th and $i_2$-th components, respectively, with
$i_1\neq i_2$. Abusing the notation, we will denote by the product
$w_1w_2\cdots w_n$ the original link $L$, by $w_i$ the link obtained
from replacing the component $L_1$ by $w_i$, and by $w_iw_j$ the link
obtained from replacing the component $L_1$ by a \lq\lq nice" band-sum
of $w_i$ and $w_j$, with $i,j$ distinct, on the disk $B$ away from other small
circles. See Figure 3.1.

\begin{figure}
\caption{\protect\footnotesize\protect\sf
The standard picture of a component in an ASL.}
\end{figure}

\bigskip

\begin{lemma} In ${\cal V}_{k-1}/{\cal V}_k$, we have
\begin{equation}
 w_1w_2\cdots w_n = \sum_{1\leq
i<j\leq n} w_iw_j-(n-2)\sum_{1\leq i\leq n} w_i.
\end{equation}
\end{lemma}

\begin{figure}
\caption{\protect\footnotesize\protect\sf
The basic relations in $\cal V_{k-1}/\cal V_k$.}
\end{figure}

\bigskip

\begin{proof} In Figures 3.2 and 3.3, we see local pictures of
$+1$-framed ASL's with
$\leq k$ components. Two basic relations in $\cal V_{k-1}/\cal V_k$, which are
consequences of Lemma 2.2, are shown in Figure 3.2. The
lemma is then proved inductively following Figure 3.3. There are three
differences in the last expression in Figure 3.3. The first difference is,
using
the second basic relation $n-1$ times,
$$\sum_{1\leq i\leq{n-1}}w_iw_n-\sum_{1\leq i\leq{n-1}}w_i-(n-1)w_n.$$
The second difference is the ASL obtained by replacing $L_1$ with $w_1\cdots
w_{n-1}$, which equals to, inductively,
$$\sum_{1\leq i<\leq{n-1}}w_iw_j-(n-3)\sum_{1\leq i\leq{n-1}}w_i.$$
Finally, the third difference is simply $w_n$. Put all these three differences
together, we get the conclusion of Lemma 3.1.
\end{proof}

By Lemma 3.1, we only need to understand to what extent the component $L_1$ in
each link $w_i$ or $w_iw_j$, $i\neq j$, is \lq\lq localized".

For the link $w_i$, there is only one triple $\mu$-invariant involving
the index 1: $\mu(1i_1i_2)=\pm1$. Locally, it can be realized as in
Figure 3.5 (a).

For the link $w_iw_j$, $i\neq j$, we need to consider several cases separately.

\medskip
\noindent{\it Case 1:} No non-zero triple $\mu$-invariant involving
the index 1. In this case, $w_iw_j\equiv1$ modulo commutators of
weight $\geq3$. So $w_iw_j=0$ in ${\cal V}_{k-1}/{\cal V}_k$.

\medskip
\noindent{\it Case 2:} Only one non-zero triple $\mu$-invariant
involving the index 1: $\mu(1i_1i_2)=\pm2$. Locally, it can be
realized as in Figure 3.5 (b).

\medskip
\noindent{\it Case 3:} Two non-zero triple $\mu$-invariant involving
the index 1: $\mu(1i_1i_2)\mu(1i_2i_3)=\pm1$ with $i_1,i_2,i_3$
distinct. Locally, they can be realized as in Figure 3.5 (c) and (d).

\medskip
\noindent{\it Case 4:} Two non-zero triple $\mu$-invariant involving
the index 1: $\mu(1i_1i_2)\mu(1i_3i_4)=\pm1$ with $i_1,i_2,i_3,i_4$
distinct. Locally, they can be realized as in Figure 3.5 (e) and (f).

For the link in Figure 3.5 (a), we may assume that $\mu(1i_1i_2)=1$
since it doesn't matter how the first component is oriented in ${\cal
V}_{k-1}/{\cal V}_k$. Similarly, in Figure 3.5 (c) and (e), we may assume
that the triple $\mu$-invariants involving the index 1 all equal to
1. Finally, as proved in proved in [O2], it is not hard to see that links in
Figure 3.5 (b), (d) and (f) are constant multiples of links in Figure 3.5 (a),
(d) and (e), respectively, in ${\cal V}_{k-1}/{\cal V}_k$.
Thus we have achieved those
two properties in Theorem 2.2 for triple $\mu$-invariant involving the
index 1.

Once we have localized the first component, we proceed to localize
the second component. If you look at the process of localizing a component,
you will find out that the localization of the second component
will effect $\mu(1i_1i_2)$ involving the index 2 only by possibly
change a non-zero one to zero. Keep localizing each component and this
completes the proof of Theorem 2.2.

\begin{figure}
\caption{\protect\footnotesize\protect\sf
 Prove Lemma 3.1 inductively.}
\end{figure}
\bigskip

To simplify the exposition, by {\it simple graphs}, we will mean
graphs satisfy (1) the valence at each vertex is either 1 or 3; and
(2) every edge connecting two distinct vertices.

\section{Induced knot invariants: The proof of Theorem 1.1 and 1.2}

Let $\lambda$ be a \ZHS invariant. Let $L$ be an ASL in $S^3$ with
$+1$ framings. We define
\begin{equation}
\psi_\lambda(L)=\sum_{ L' \subset  L}(-1)^{\#L'}
\lambda(S^3_{L'}).
\end{equation}

\begin{lemma} A \ZHS invariant $\lambda$ is of order $\leq k$ iff
$\psi_{\lambda}(L)=0$ for every ASL $L$ in $S^3$ with $\#L>k$.
\end{lemma}

\begin{proof} This comes from the definition (1.1) and Lemma 2.1.
\end{proof}

Analogous to the theory of finite type knot invariants, we may think
of ASL's in $S^3$ as \lq\lq singular $\Bbb Z$-homology 3-spheres" and
simple graphs as \lq\lq chord diagrams". For a finite type \ZHS
invariant $\lambda$, $\psi_\lambda$ may then be thought of \lq\lq
derivatives" of $\lambda$. The following two lemmas have their counterparts in
the theory of finite type knot invariants.

\addtocounter{section}{-1}
\addtocounter{figure}{3}
\begin{figure}
\caption{\protect\footnotesize\protect\sf
 Localization of a component.}
\end{figure}
\bigskip
\addtocounter{section}{1}
\addtocounter{figure}{-4}

\begin{lemma} Let $\lambda$ be a \ZHS invariant of order $\leq k$. Then
$\lambda$ is
determined by $\psi_{\lambda}$ on ASL's in $S^3$ with $\leq k$
components.
\end{lemma}

\begin{proof} Let $M=S^3_L$ for a $\pm1$-framed ASL $L$ in $S^3$. We prove the
lemma inductively on the number of components of $L$.

If $L=\emptyset$, then $M=S^3=(S^3,\emptyset)$. So the lemma holds. Now assume
that $\#L=n$. Then,
$$\lambda(S^3_L)=(-1)^{n-1}\left(\psi_\lambda(L)-\sum_{L'\subset
L,L'\neq L}\lambda(S^3_{L'})\right).
$$
Inductively, $\lambda(S^3_{L'})$ for $L'\subset L,\,L'\neq L$ is determined by
$\psi_\lambda$ on ASL's in  $S^3$ with $\leq k$ components. Since
$\psi_\lambda(L)=0$ if $n>k$, we see that $\lambda(S^3_L)$ is also determined
by $\psi_\lambda$ on ASL's in $S^3$ with $\leq k$ components.
\end{proof}

\begin{lemma}
Let $L$ be an ASL in $S^3$ and $O$ be an unknot separated from $L$, then
$\psi_\lambda(L \cup O) =0$.
\end{lemma}

\begin{proof} We have
$$
\begin{aligned}
\psi_\lambda(L\cup O)& = \sum_{L' \subset L, O
        \notin L'} (-1)^{\#L'}\lambda(S^3_{L'}) +
        \sum_{L' \subset L , O
        \in L'} (-1)^{\# L'}\lambda(S^3_{L'})\\
&=\sum_{L' \subset L , O
        \notin L'} (-1)^{\#L'}\lambda(S^3_{L'})+
\sum_{L' \subset L , O
        \notin L'} (-1)^{\# L'+1}\lambda(S^3_{L'\cup O}).
\end{aligned}$$
Notice that $S^3_{L'\cup O} = S^3_{L'}$ for any sublink $L'$ of $L$ such that
$O\notin L'$. Hence we have
$$ \psi_\lambda(L \cup O) = \sum_{L' \subset  L ,  O
        \notin  L'} (-1)^{\# L'}\lambda(S^3_{L'}) -
        \sum_{ L' \subset L ,O\notin L'}
        (-1)^{ L'}\lambda(S^3_{L'}) =0.$$

This finishes the proof.
\end{proof}

Restricting to the case $\#L=1$, $\psi_\lambda$ becomes a knot
invariant. We emphasize this situation by having a separate notation:
For a knot $K$ in $S^3$, frame it by $+1$ and denote
\begin{equation}
v_\lambda(K)=\lambda(S^3)-\lambda(S^3_K).
\end{equation}

Notice that this definition of $v_\lambda$ differs slightly with the
definition, which is given in (1.2), used in the statement of Theorem 1.1.
Obviously, the conclusions of Theorem 1.1 under these two definitions are
equivalent. So we may use the current definition in the proof of
Theorem 1.1.

\begin{lemma} If $\lambda$ is a \ZHS invariant of order $\leq k$, $v_\lambda$
is a
knot invariant of order $\leq k-1$.
\end{lemma}

\begin{proof}
This follows from Ohtsuki's first local move in Figure 2.1. For the
evaluation of $v_\lambda$ on a knot with $k$ double points is equal to
$\psi_\lambda(L)$ for a $\pm1$-framed ASL $L$ in $S^3$ with $k+1$
components, which is zero by definition of $\lambda$.
\end{proof}

Given any chord diagram $C$ with $k-1$ chord, one may obtain an ASL
with $k$ components in $S^3$ in the way depicted in Figure 4.1. We
will first replace each chord by a thin untwisted band intersecting in
its interior the outer circle at the end points of the original
chord. These bands are disjoint from each other. We then take the
boundaries of these bands together with the outer circle to form a ASL
with $k$ components. It is a well-defined element in ${\cal
V}_{k-1}/{\cal V}_k$ by Ohtsuki's second local move. We denote by $L(C)$
this ASL associated with a chord diagram $C$. The following lemma was
observed by Garoufalidis [G].

\begin{figure}
\caption{\protect\footnotesize\protect\sf
 A chord diagram, its associated ASL and a singular
knot respecting it.}
\end{figure}

\bigskip

\begin{lemma} The weight system for $v_\lambda$ is given by
$$v_\lambda(C) = \psi_\lambda(L(C))$$
for every chord diagram $c$ with $k-1$ chords.
\end{lemma}

\begin{proof} Choose a ``good'' immersion for $C$ (which is a singular knot
respecting $C$). This is to be done in the  following way
[B-L]. Pick a point on the outer circle and follow it
counterclockwise. When you come to a chord, follow the chord till it
meets the other edge of the outer circle, fix this as a double point,
loop back under the outer circle and follow the chord back to where it
began, continue around the outer circle, till you have returned to the
starting point. See Figure 4.1. Notice now that the immersion can always
be chosen in the following way. When each double point is resolved
positively, a loop is created that can be isotoped back to the outer
circle. The net effect being as if that chord in the chord diagram
was erased. Further, when each double point is resolved negatively,
the twist created in the outer circle is a negative one. The same
effect can be certainly achieved by either performing or not
performing surgery on the component of $L(C)$ corresponding to that
chord. Thus,
$$\begin{aligned}
v_\lambda(C)&=\sum_{L'\subset L(C)-L_0}(-1)^{\# L'}\lambda(S^3_{L'})-
\sum_{L'\subset L(C)-L_0}(-1)^{\#L'}\lambda(S^3_{L'\cup L_0})\\
&=\psi_\lambda(L(C)).
\end{aligned}
$$
\end{proof}

We now use Lemma 3.1 to calculate the weight system for
$v_\lambda$. Notice that the outer circle, call it $L_0$, in the ASL
$L(C)$ associated with a chord diagram $C$ is quite special. If we
delete $L_0$, $L(C)$ will fall apart into an unlink. Therefore, $L_0$
represents an element in the free group $\pi_1(S^3\setminus(L(C)-L_0))$
with $k-1$ generators. Moreover, we can write out this element as a
word in generators easily as follows. First,
arbitrarily number the chords by $1,\dots,k-1$, and orient them. For
each chord $i$, place the symbol $i$ on the outgoing endpoint of the
chord and $\bar i$ on the incoming endpoint. Pick a base point on the
outer circle and read off the symbols encountered when traversing the
circle counterclockwise. This is the word we are looking for. This
word can be expressed as a product of simple commutators of the
meridians of $L(C)-L_0$ modulo commutators of weight $\geq3$. If the
number of simple commutators is $>2$, we may use Lemma 3.1 to
decompose this product of simple commutators into a sum of products of
at most two simple commutators.  Since there is no non-zero triple
$\mu$-invariant other than these involving the index $0$ (of the
component $L-L_0$, such a decomposition suffices to express $L(C)$ as a
sum of simple graphs.

As an example, we first reprove Garoufalidis' result in the following lemma.

\begin{figure}
\caption{\protect\footnotesize\protect\sf
 The bubble and the switch.}
\end{figure}

\bigskip

\begin{lemma} ${\cal O}_4={\cal O}_3$.
\end{lemma}

\begin{proof} There is only one simple graph of 4 edges, ``the bubble" (see
Figure
4.2),
which is possibly non-zero in ${\cal V}_3/{\cal V}_4$. Suppose there is a
\ZHS invariant $\lambda$ of order $\leq4$ and let us consider the
induced knot invariant of order $\leq3$.  There are two chord diagrams
$C_1$ and $C_2$ with 3 chords which may have non-zero weights for an
invariant of order 3. There are shown in Figure 4.3. The 4-term
relation in the theory of finite type knot invariants implies that
$v_\lambda(C_2)=2v_\lambda(C_1)$.  Thus, we must have
$$\lambda(L(C_2))=2\lambda(L(C_1)).$$

Now consider their respective ASL's. It
is not to hard to see that $L(C_1)$ can be decomposed directly as a
simple graph, the bubble.  Let us now decompose $L(C_2)$ into simple
graphs. We label the chords of $C_2$ by 1,2,3 and obtain the following
word
$$123\bar 1\bar 2\bar 3 = [1,23][2,3] \sim [1,2][1,3][2,3]$$
By Lemma 3.1, we can decompose this word as
$$[1,2][1,3]+[1,2][2,3]+[1,3][2,3]-[1,2]-[1,3]-[2,3].$$ The simple
graphs corresponding to the last three words are zero in ${\cal V}_3/{\cal
V}_4$ since they all have an isolated edge. The simple graphs
corresponding to the first three words are all equal to the bubble. Thus
$$L(C_2)=3L(C_1)$$
in ${\cal V}_3/{\cal V}_4$. This is possible only when
$\lambda=0$ on the bubble. Therefore, $\lambda$ is of order $\leq3$.
\end{proof}

\begin{figure}[h]
\caption{\protect\footnotesize\protect\sf
 Chord diagrams with 3 chords.}
\end{figure}

\begin{ppf} For $k=4$, the proof of Theorem 1.1 is given by Lemma 4.6.
We jump to the case $k>5$ first. Let $\lambda\in {\cal O}_k$ for
$k>5$. Let us evaluate the weight system of $v_\lambda$. For each
chord diagram $C$ with $k-1>4$ chords, decompose the outer circle
$L_0$ in $L(C)$ using Lemma 3.1 into a sum of products of at most two
simple commutators. Since there are more than four chords, there will
be at least one isolated edge in each of the simple graphs in the
decomposition coming from Lemma 3.1 and
$v_\lambda(C)=\lambda(L(C))=0$.

Finally, we come to the proof of Theorem 1.1 in the case when
$k=5$. Again, in this case, there is only one simple graph with 5
edges, ``the switch'', which is possibly non-zero in ${\cal V}_4/{\cal
V}_5$. See Figure 4.2. Unlike in the proof of Lemma 4.6, the induced weight
system on chord diagram with 4 chords is actually consistent as a
one-dimensional
vector space, whose basis element is the switch.

Nevertheless, as we will see, the local moves in Figure 2.1 put more
constraints on induced knot invariants than merely the 4-term relation. Let
us consider the following chord diagram $C$ with 4 chords. (Figure 4.4)
Its word is given by
$$1234\bar 1\bar 2\bar 3\bar 4 = [1,234][2,34][3,4] \sim [1,2][1,3][1,4][2,3]
[2,4][3,4]$$
Using Lemma 3.1, this evaluates to three times the switch. On the other hand,
if, instead of the outer loop, we choose to apply Lemma 3.1 to one of the loops
associated with a chord, we will arrive at the following diagrammatic
expansion. See Figure 4.4. We first change the position of the loop associated
with the chord marked by 2 so that it is in a position for us to apply Lemma
3.1.
We get the next linear combination of ASL's by applying Lemma 3.1. Notice the
last
two ASL's in this linear combination are associated with, respectively, chord
diagrams
$1234\bar1\bar3\bar2\bar4$ and $12\bar134\bar2\bar3\bar4$. They can
easily be evaluated to 2 times the switch and the switch respectively. The
first ASL, though, is not associated with any chord diagram. But it still has
the property that when we drop the outer circle, all the other components fall
apart. Thus, the same method we used to evaluate a chord diagram can be used
to evaluate this ASL. We use the first diagram in the last row in Figure 4.4
to
denote this ASL. The small circle at the intersection of chord 2 and chord 4
indicates that when we want to recover the original ASL, the band associated
with chord 2 runs through the band associated with chord 4. Hence we may read
off the word represented by the outer circle as
$$1234\bar 1\bar 2\bar 32\bar 4 \bar2 = [1,234][2,34][3,4][4,2]$$
and it evaluates to 4 times the switch. Thus we get the equality that 3
times the switch is equal to 5 times the switch. This is possible only when any
\ZHS invariant of order $\leq 5$ evaluates zero on the switch. We therefore
have
proved Theorem 1.1 in the case when
$k=5$ as well as Theorem 1.2.
\end{ppf}

\begin{figure}
\caption{\protect\footnotesize\protect\sf
 The proof of Theorem 1.1 when $k=5$ as well as
Theorem
1.2.}
\end{figure}

\bigskip

The proof of Theorem 1.1 in the case when $k>5$ suggests that knots
are too simple to hold up all finite type \ZHS invariants. On the
other hand, Lemma 4.6 and Theorem 1.2 tell us that a detailed
knowledge about induced link invariants is of great help in
understanding the original 3-manifold invariants. A natural question
will be: What properties characterizes the induced invariant on ASL's
with two components?

\end{document}